\documentclass{article}
\usepackage{amsmath,amsfonts,amssymb,amsthm}
\usepackage{graphicx} 
\usepackage{subcaption}
\usepackage[english,greek]{babel}
\usepackage{tikz}
\usepackage{url}

\usepackage{comment}

\newcommand{\bF}{\mathbb{F}}

\newcommand{\bR}{\mathbb{R}}

\newcommand{\fS}{\mathfrak{S}}

\newcommand{\interior}{\mathbf{int}}
\newcommand{\ci}{\mathbf{ci}}
\newcommand{\ic}{\mathbf{ic}}
\newcommand{\closure}{\mathbf{cl}}

\newcommand{\gdc}{\mathfrak{c}}

\newcommand{\mb}[1]{\mathbf{#1}}
\newcommand{\mr}[1]{\mathrm{#1}}

\newcommand{\sto}[1]{\widetilde{#1}}
\newcommand{\stok}[2]{\widetilde{#1}^{(#2)}}
\newcommand{\unitv}[1]{\frac{\overrightarrow{#1}}{\|\overrightarrow{#1}\|}}

\theoremstyle{definition}
\newtheorem{defn}{Definition}[section]
\newtheorem{lem}[defn]{Lemma}

\newtheorem{thm}[defn]{Theorem}
\newtheorem{cor}[defn]{Corollary}
\theoremstyle{remark}

\newtheorem{ex}[defn]{Example}

\title{
How can a cake be cut into two equal pieces?
}
\author{Byungchang So\\ \small
\url{sinwall@snu.ac.kr} \\ \small
College of Humanities, \\ \small
Seoul National University \\ \small
}
\date{}
\begin{document}
\selectlanguage{english}
\usetikzlibrary{math}
\maketitle
\begin{abstract}
In coordinate geometry, geometric objects such as lines, disks, and cubes are not perfectly evenly divided because the points on the border should belong to any one of the parts.
To mathematically formalize such division, which is physically straightforward, this article proposes a novel representation of geometric objects, called ``rasped representation.''
\end{abstract}


\section{Introduction}
\begin{figure}[h]
\centering
\includegraphics{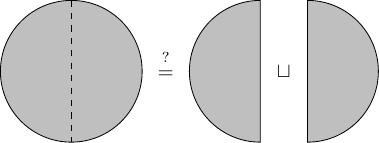}
\caption{Is a round cake cut into two equal pieces?}
\label{fig:cake}
\end{figure}

The answer is simple:
By cutting along any diameter (no need to see Figure \ref{fig:cake}).

But can the same answer be given in coordinate geometry?
Is disk $\{(x,y)\in \bR^2 : x^2+y^2 \leq r^2 \}$ equal to the \textit{disjoint union} of two congruent half-disks?
No, it is not!
If each half-disk includes its boundary, then the two half-disks are not disjoint.
If each half-disk does not include its boundary, then the two half-disks do not fill the entire disk.\footnote{%
Essentially the same question was raised by Constantin Carath\'{e}ory in \cite[p.12]{Caratheodory1963}.
In his example, a parallelogram, rather than a disk, is cut along one of its diagonals.}
A disk can be divided into two half-disks with \textit{equal areas}, but not as a disjoint union.
\begin{figure}
\centering
\includegraphics{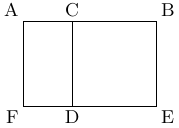}
\caption{\label{fig:elem2-3}%
A diagram for \textit{Elements} II.3.
}
\end{figure}

How can cutting a cake be described in mathematics?
This is more than a recreational question and is indeed what readers of classical Greek geometry encounter, e.g. in Euclid's \textit{Elements}, Archimedes' \textit{On Spirals}, Apollonius' \textit{Conics}.
For example, in \textit{Elements} II.3 (see Figure \ref{fig:elem2-3}),
\begin{center}
$\square\mr{ABEF}$ is equal to $\square\mr{ACDF}$ and $\square\mr{CBED}$.\footnote{%
This is a liberal translation by the author.
A more literal translation by Thomas Little Heath in \cite{HeathElements1} reads ``$\mr{AE}$ is equal to $\mr{AD}$, $\mr{CE}$.'' } 
\end{center} 
We would mathematically formulate the above statement as $\text{area}(\square\mr{ABEF}) = \text{area}(\square\mr{ACDF}) + \text{area}(\square\mr{CBED})$ or $\square\mr{ABEF} = \square\mr{ACDF} \cup \square\mr{CBED} - \mr{CD}$. 
But why not avoid the quantification by area and subtracting $\mr{CD}$, neither of which present in the original text, and instead try $\square\mr{ABEZ} = \square\mr{ACDF} \, \sqcup \, \square\mr{CBED}$ as a modern mathematical formulation?\footnote{\label{fn:equal}%
The meaning of ``equal(\textgreek{ἴσος})'' must be considered further, though.
For example, two scissor-congruent polygons are also equal'' by \textit{Common Notions} 1--3 in \textit{Elements} I: \begin{enumerate} 
\item[CN1.] Things equal to the same thing are equal to each other. 
\item[CN2.] If equal things are added to equal things, then the wholes are equal. 
\item[CN3.] If equal things are removed from equal things, then the remainders are equal. \end{enumerate} 
Interested readers may refer to \cite[\S3]{DeRisi2021} for a discussion on the concept of equality in \textit{Elements}. 
See also \cite[\S2--3]{GrattanGuinness1996} for more on differences between the mathematical style of \textit{Elements} and its modern interpretation.
}

To achieve the ``seamless division'', we do not represent each position in a space by a single element but allow it to be divisible further. 
In the cake example, each position on the cut line will have a left part and a right part, each occupied by a half-disk. 
Consequently, geometric objects are represented not as sets of points\footnote{%
David Hilbert's \textit{Grundlagen} \cite{Grundlagen}, a central work in formalization of Euclidean geometry, already assumes lines and polygons to be sets of points. 
Interested readers can find a more recent exposition of \textit{Grundlagen} in \cite{Hartshorne2000}.
} 
but as sets of something finer. 
The new representation, which we will call ``rasped'' representation (see Section \ref{sec:motiv} for an explanation of this term), possesses the following properties, given coordinate space $\bR^d$.
\begin{itemize}
\item Geometric objects are represented as subsets of $\bR^d \times S_{d,k}$ (with a temporary notation), where $S_{d,k}$ varies with the dimension $k$ of the given geometric object ($k=0$ for point, $k=1$ for line, $k=2$ for polygon, etc.) as well.
See Definition \ref{def:sto-basic} and subsequent formula \eqref{eq:recipe} for details. 

\item Rasped representation is defined inductively. 
Roughly speaking, the rasped representation of $X$ binds to each point $x$ of $X \subseteq \bR^d$ the rasped representation of $U_x X$, the infinitesimally local shape at $x$.  
See equations \eqref{eq:recipe} and \eqref{eq:rasp-curv} for details. 

\item Elementary geometric objects such as polygons, sectors, polyhedra, cones, etc. admit rasped representations. 
For these objects, the rasped representation is, in a sense, equivalent to conventional representation in coordinate geometry.  
See Theorem \ref{thm:rasp-curv} for details. 
\end{itemize}
Hence the (genuine) answer to the title is: Cut the rasped representation, which is, quite literally, a piece of cake!

\paragraph{Organization}{%
The contents of each section are as follows.
\begin{itemize}
\item Section 2 illustrates rasped representation through heuristic arguments.
\item Section 3 examines several examples of rasped representation.
\item Section 4 provides rigorous accounts of rasped representations and seamless division.
\item Section 5 closes this article.
\end{itemize}
Sections 2 and 3 assume familiarity with basic linear algebra and calculus.
Section 4 uses concepts in Riemannian geometry and set-theoretic topology.
}

\paragraph{Notations}{%
The following symbols are used frequently in this article.
\begin{itemize}
\item The alphabets in non-boldface roman ($\mr{A}$, $\mr{B}$, $\mr{C}$, $\cdots$) are reserved for the names of points.

\item $\mb{A}$, $\mb{V}$, and $\mb{S}$ are: $\mb{A} = \mb{V} = \mb{\bR}^d$, $\mb{S} = \{(x_1,\cdots,x_d)\in \bR^d: x_1^2 + \cdots + x_d^2 = 1\}$ where $d \in \mathbb{N}$.
The dot product of $a,b \in \mb{V}$ is denoted by $\langle a,b\rangle = a_1b_1 + \cdots + a_d b_d$.

The content of this article remains valid even if these symbols are understood as follows:
$\mb{V}$ is an $\bF$-vector space equipped with inner product $\langle \cdot,\cdot\rangle$ and its induced norm $| \cdot |$.
$\mb{A}$ is an $\bF$-affine space associated with $\mb{V}$.
$\mb{S} := \{ v \in \mb{V} : |v | = 1\}$ is the unit sphere in $\mb{V}$.

\item $[a,b]$ is the closed interval with endpoints $a$, $b$ in $\bR$ or unit circle $\{(x,y)\in\bR^2: x^2+y^2=1\}$.
In the latter case, this refers to the trace of the counterclockwise trip from $a$ to $b$.

\item $\sqcup$ and $\bigsqcup$ refer to disjoint union of sets.
For example, $X = Y \sqcup Z$ is equivalent to $X = Y \cup Z$ and $Y \cap Z = \varnothing$.

\item $\interior\,X$ and $\closure\,X$ refer to the interior and closure, respectively, of a subset $X$ of a topological space.
Iterated operations are abbreviated as $\ci (X) := \closure(\interior\, X)$, $\ic (X) := \interior(\closure\, X)$, and so forth.
\end{itemize}
}

\section{A rough sketch} \label{sec:motiv}

\textit{%
This section proceeds with heuristic arguments that may not be mathematically rigorous.
}

\begin{figure}
\centering
\begin{subfigure}{0.42\textwidth}
\centering
\includegraphics{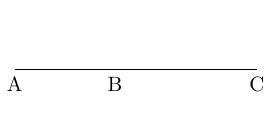}
\caption{$\sto{\mr{AB}} \sqcup\sto{\mr{BC}} = \sto{\mr{AC}}$}
\label{subfig:seamless-1d}
\end{subfigure}
\begin{subfigure}{0.42\textwidth}
\centering
\includegraphics{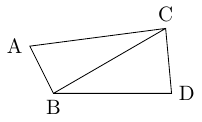}
\caption{$\sto{\triangle\mr{ABC}} \sqcup \sto{\triangle\mr{BCD}} = \sto{\square\mr{ABCD}}$}
\label{subfig:seamless-2d}
\end{subfigure}
\caption{Examples of seamless division.}
\end{figure}

The goal of this article is to represent elementary geometric objects (e.g., lines, rectangles, cubes) in a way that enables ``seamless division'' (rigorous definition of this term will appear later in Section \ref{sec:rasp-curv}) as follows. 
In the case of lines, if two line segments $\mr{AB}$ and $\mr{BC}$ meet at the point $\mr{B}$ only, then the sets $\sto{\mr{AB}}$, $\sto{\mr{BC}}$, and $\sto{\mr{AC}}$ representing each line segment should satisfy $\sto{\mr{AB}} \sqcup\sto{\mr{BC}} = \sto{\mr{AC}}$ (see Figure \ref{subfig:seamless-1d}).
In the case of planar regions, if two triangles $\mr{ABC}$ and $\mr{BCD}$ meet each other along the line segment $\mr{BC}$ only, then the sets $\sto{\triangle \mr{ABC}}$, $\sto{\triangle\mr{BCD}}$, and $\sto{\square\mr{ABCD}}$ representing each triangle or quadrilateral should satisfy $\sto{\triangle\mr{ABC}} \sqcup \sto{\triangle\mr{BCD}} = \sto{\square\mr{ABCD}}$ (see Figure \ref{subfig:seamless-2d}). 

We will call the representation \textit{rasped} for the following reason.
For example, if the triangles $\mr{ABC}$ and $\mr{BCD}$ above are represented as closed subsets of $\bR^2$, they overlap along the line segment $\mr{BC}$. 
For seamless division to hold, their boundaries should be slightly rasped (shaved off), but not excessively so that $\mr{BC}$ remains in the union of the two triangles. 

Let us examine basic examples.
In what follows, symbols with tilde(``$\sim$'') always refer to a rasped representation of an object.

\begin{figure}
\centering
\begin{subfigure}{0.42\textwidth}
\centering
\includegraphics{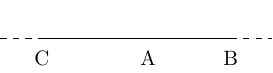}
\caption{\label{subfig:1da}}
\end{subfigure}
\begin{subfigure}{0.42\textwidth}
\centering
\includegraphics{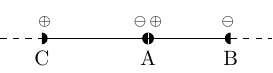}
\caption{\label{subfig:1db}}
\end{subfigure}
\caption{Rasped representation of line segments.}
\label{fig:motiv-1d}
\end{figure}

\paragraph{Step 1: Line segments}{%
Suppose first that three points $\mr{C}$, $\mr{A}$, and $\mr{B}$, from left to right, lie on a infinite straight line $\mb{A}$ (see Figure \ref{subfig:1da}). 
Let $\sto{\mr{CA}}$, $\sto{\mr{AB}}$, and $\sto{\mr{CB}}$ be the rasped representations of line segments $\mr{CA}$, $\mr{AB}$, and $\mr{CB}$, respectively, and $\sto{\mb{A}}$ the ambient set containing these three representations. 
Then seamless division
\begin{equation*}
\sto{\mr{CA}} \ \sqcup\  \sto{\mr{AB}} \ =\  \sto{\mr{CB}}
\end{equation*}
must hold. 
Since the point $\mr{A}$ lies on the segment $\mr{CB}$, which is now represented as $\sto{\mr{CA}} \sqcup \sto{\mr{AB}}$, the position of point $\mr{A}$ should be divided into two parts, say $\ominus$ and $\oplus$, where $\ominus$ is occupied by $\sto{\mr{CA}}$ and $\oplus$ by $\sto{\mr{AB}}$ (see Figure \ref{subfig:1db}). 
Likewise, the position of point $\mr{B}$ should also have two parts $\ominus,\oplus$, and only $\ominus$  is occupied by $\sto{\mr{AB}}$. 
Assuming homogeneity, every other point on the segment $\mr{AB}$ also has two parts, both occupied by $\sto{\mr{AB}}$. 

Thus, we obtain the following quasi-formulas for the case $\dim \mb{A} = 1$, where the question marks mean that no rigorous definitions are provided yet for the equalities to make sense.
\begin{enumerate}
\item $\sto{\mb{A}} \stackrel{?}{=} \mb{A} \times (\ominus \sqcup \oplus)$

\item $\sto{\mr{AB}} 
\stackrel{?}{=} \{\mr{A}\} \times \oplus 
\ \sqcup\  
\bigsqcup_{\mr{C} \text{ lies inside } \mr{AB}}
\{\mr{C}\} \times (\ominus \sqcup \oplus)
\ \sqcup\  \{\mr{B}\} \times \ominus$
\end{enumerate}
Note that a similar conclusion can be drawn if $\mb{A}$ and line segments are replaced with a circle and arcs.

\begin{figure}
\centering
\begin{subfigure}{0.32\textwidth}
\centering
\includegraphics{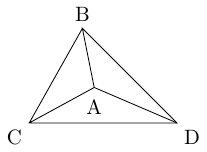}
\caption{\label{subfig:2da}}
\end{subfigure}
\begin{subfigure}{0.32\textwidth}
\centering
\includegraphics{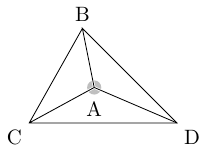}
\caption{\label{subfig:2db}}
\end{subfigure}
\begin{subfigure}{0.32\textwidth}
\centering
\includegraphics{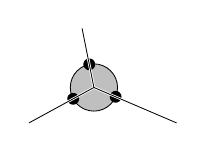}
\caption{\label{subfig:2dc}}
\end{subfigure}
\caption{Rasped representation of triangles.}
\label{fig:motiv-2d}
\end{figure} 
}

\paragraph{Step 2: Triangles}{%
Next, suppose that a point $\mr{A}$ is placed inside a triangle $\mr{BCD}$ in a plane $\mb{A}$ (see Figure \ref{subfig:2da}). 
Let $\sto{\triangle \mr{ABC}},\sto{\triangle \mr{ACD}},\sto{\triangle \mr{ADB}}$, and $\sto{\triangle \mr{BCD}}$ be the rasped representations of triangles $\mr{ABC}$, $\mr{ACD}$, $\mr{ADB}$, and $\mr{BCD}$, respectively. 
Then seamless division
\begin{equation*}
\sto{\triangle \mr{ABC}} \ \sqcup\ 
\sto{\triangle \mr{ACD}} \ \sqcup\ 
\sto{\triangle \mr{ADB}}\ =\  \sto{\triangle \mr{BCD}}
\end{equation*}
must hold. 
Since the point $\mr{A}$ lies on the triangle $\mr{BCD}$, which is now represented as $\sto{\triangle \mr{ABC}} \sqcup \sto{\triangle \mr{ACD}} \sqcup \sto{\triangle \mr{ADB}}$, the position of the point $\mr{A}$ should be divided into at least three parts occupied by each of three triangles $\mr{ABC},\mr{ACD}$ and $\mr{ADB}$ (see Figure \ref{subfig:2db}). 

Taking all possible positions of $\mr{B}$, $\mr{C}$, and $\mr{D}$ into account, the position of $\mr{A}$ is likely to be divisible into infinitely many parts, one for each element of the circle $S_{r}(\mr{A}) := \{x\in\bR^2: \|x-\mr{A}\|=r\}$. 
If this is the case, the parts occupied by the triangles $\mr{ABC}$, $\mr{ACD}$, and $\mr{ADB}$ would correspond to three arcs of $S_{r}(\mr{A})$ with central angles $\angle\mr{BAC}$, $\angle\mr{CAD}$, and $\angle\mr{DAB}$, respectively. 
Since the three triangles are mutually disjoint, the three occupied parts should also be mutually disjoint. 
It is as though the entire circle were equal to the disjoint union of the three arcs -- the current situation reminds us of the Step 1 above! 
Then, the circle $S_{r}(\mr{A})$ should be replaced by its rasped representation $\sto{S_{r}(\mr{A})}$ (see Figure \ref{subfig:2dc}). 
Accordingly, the parts of the position of $\mr{A}$ occupied by the triangle $\mr{ABC}$ correspond to the rasped representation of the arc with central angle $\angle\mr{BAC}$. 
A similar argument also applies to other points lying in the triangle $\mr{ABC}$. 

Thus, we obtain the following quasi-formula for the case $\dim \mb{A} = 2$.
\begin{enumerate}
\item $\sto{\mb{A}} \stackrel{?}{=} \mb{A} \times \sto{\mb{S}}$, where $\sto{\mb{S}}$ is the rasped representation of a circle $\mb{S}$. 

\item $\sto{\triangle \mr{ABC}} \stackrel{?}{=} \{\mr{A}\} \times \sto{(\text{arc for }\angle\mr{BAC})} 
\,\sqcup\, \{\mr{B}\} \times \sto{(\text{arc for }\angle\mr{CBA})}$ \\
$\phantom{\sto{\triangle \mr{ABC}} = {}} 
\sqcup\, \{\mr{C}\} \times \sto{(\text{arc for }\angle\mr{BCA})}\sqcup \cdots$
\end{enumerate}
Note that a similar conclusion can be drawn for other polygons and even if $\mb{A}$ and polygons are replaced with $2$-dimensional sphere and geodesic polygons. 

\paragraph{Step 0: Point}{%
Meanwhile, the situation is trivial for the $0$-dimensional affine space (one-point set) and the sphere (two-point set), because seamless division already holds. 
We may just let $\sto{\mb{A}} = \mb{A}$ and $\sto{\mb{S}} = \mb{S}$.
}

\paragraph{Generalization}{%
Based on Steps 0--2, what can we expect about rasped representation in affine spaces and spheres of arbitrary dimension?

The division of individual positions in an affine space or a sphere of dimension $d$ would be represented by a sphere of dimension $d-1$. 
This leads to, for an affine space $\mb{A}$
\begin{equation}\label{eq:rasp-affine-quasi}
\sto{\mb{A}} \stackrel{?}{=} \mb{A} \times \sto{\mb{S}},
\end{equation}
where $\mb{S}$ is a sphere in the associated vector space of $\mb{A}$, and for a sphere $\mb{S}'$ in a vector space $\mb{V}'$
\begin{equation}\label{eq:rasp-sphere-quasi}
\begin{aligned}
\sto{\mb{S}'} 
& \stackrel{?}{=} \bigsqcup_{\mr{A} \in \mb{S}'} \{\mr{A}\} \times \sto{(\text{sphere of the tangent space at } \mr{A})} \\
& \stackrel{?}{=} \bigsqcup_{\substack{v \in \mb{V}'\\ \|v \|=1}} \{v\} \times \sto{(\text{sphere of } \{v\}^\perp)}.
\end{aligned}
\end{equation} 
Equations \eqref{eq:rasp-affine-quasi} and \eqref{eq:rasp-sphere-quasi} form a recurrence relation with initial condition given in Step 0 above.
We are now ready to establish a rigorous definition of the ambient spaces of rasped representation. 

\begin{defn}\label{def:sto-basic}
Let $\mb{A}$ be an affine space, $\mb{V}$ its associated vector space, and $\mb{S}$ the unit sphere in $\mb{V}$. 
In addition, let $k$ be an integer such that $0 \leq k \leq \dim\mb{A}$. 
We define the sets $\stok{\mb{S}}{k-1}$ $(k \geq 1)$ and $\stok{\mb{A}}{k}$ as
\begin{equation*}
\stok{\mb{S}}{k-1} = 
\{(v_0,\cdots,v_{k-1}) \in \mb{V}^{k}: v_0,\cdots,v_{k-1} \text{ are orthonormal} \}
\end{equation*}
and
\begin{equation*}
\stok{\mb{A}}{k} 
=
\begin{cases}
\mb{A}, & k = 0 \\
\mb{A} \times \stok{\mb{S}}{k-1}, & k > 0. 
\end{cases}
\end{equation*}
Elements of either sets will be denoted e.g. by $(v_0 ; v_1 , \cdots, v_{k-1}) \in \stok{\mb{S}}{k-1}$ and $(x; v_0,\cdots, v_{k-1}) \in \stok{\mb{A}}{k}$ to emphasize the first entry.
\end{defn}
The spaces $\sto{\mb{A}}$ and $\sto{\mb{S}}$ appeared earlier in this section are actually $\sto{\mb{A}} = \stok{\mb{A}}{d}$ and $\sto{\mb{S}} = \stok{\mb{S}}{d-1}$ with $d=\dim{\mb{A}}=\dim{\mb{S}}+1$. 
The sets $\stok{\mb{A}}{k}$ with $k < \dim \mb{A}$ and $\stok{\mb{S}}{k-1}$ with $k \le \dim \mb{S}$ will turn out to be ambient sets for lower-dimensional objects.
For example, under rasped representation, polygons of an infinite plane $\mb{B}$ are subsets of $\stok{\mb{B}}{2}$, line segments of an infinite line $\mb{A} \subset\mb{B}$ are subsets of $\stok{\mb{A}}{1} \subset \stok{\mb{B}}{1}$. 

Next, the formula for rasped objects is constructed as follows. 
Let $X \subset \mb{A}$ be a closed subset (i.e. contains boundary) corresponding to a geometric object of dimension $k \leq \dim\mb{A}$. 
If $k=0$, seamless division holds without any modification as seen in Step 0, so we let $\sto{X} = X \subset \stok{\mb{A}}{0} = \mb{A}$. 
If $k \geq 1$, given $x \in X$, the part occupied by rasped representation $\sto{X} \subset \stok{\mb{A}}{k}$ corresponds to\footnote{%
The limit $\lim_{r \searrow 0}$ was unnecessary for rectilinear shapes in Step 2 but should be included for curvilinear shapes. 
The precise meaning of the expression will be clarified later in Section \ref{sec:rasp-curv}.
}
\begin{align*}
\lim_{r \searrow 0}{\frac{X-\{x\}}{r} \cap \mb{S}} 
& \stackrel{?}{:=} \lim_{r \searrow 0} \left\{ \frac{y-x}{r} \in \mb{V}:y \in X \right\} \cap \mb{S}\\
& \phantom{:}\stackrel{?}{=} \lim_{r \searrow 0}\{v \in \mb{S}: x + rv \in X \}.
\end{align*}
Recall that in Step 2, the occupation of each triangle at each position corresponded to rasped representations of arcs and not to arcs themselves. 
In the current case, this implies
\begin{equation*}
\sto{X} \cap (\{x\} \times \stok{\mb{S}}{k-1})
\stackrel{?}{=} \left[ \lim_{r \searrow 0}{\frac{X-\{x\}}{r}}\cap \mb{S} \right]^{\sim}.
\end{equation*}
Thus, the general formula should be
\begin{equation} \label{eq:recipe}
\sto{X} \stackrel{?}{=} 
\begin{cases}
X, & k = 0 \\ 
\bigsqcup _{x \in X} \{x\} \times \left[ \lim_{r \searrow 0}\frac{X-\{x\}}{r}\cap \mb{S} \right]^{\sim}, & k \geq 1,
\end{cases}
\end{equation}
and a similar formula will apply when $X$ is a subset of a sphere. 
The recipe \eqref{eq:recipe} is recursive: if $X$ is ``nice'' enough so that $\frac{X-\{x\}}{r}\cap \mb{S}$ is of dimension $k-1$ and converges to another ``nice'' object admitting rasped representation in $\stok{\mb{S}}{k-1}$, we obtain the rasped representation $\sto{X} \subset \stok{\mb{A}}{k}$. 
We will come back to the rigorous version of recipe \eqref{eq:recipe} in Section \ref{sec:rasp-curv}. 
}

\section{Examples}

Let us observe how the recipe \eqref{eq:recipe} operates and how seamless division holds in basic examples. 

In the remaining part of this article, rasped representations are found for closed subsets of either $\mb{A}$ or $\mb{S}$.
Also, we will call rasped representations of various objects rasped lines, rasped polygons, etc. 

\begin{figure}
\centering
\includegraphics{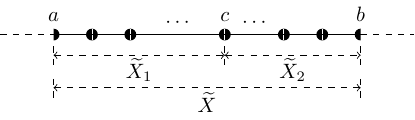}
\caption{A rasped line segment and its division.}
\label{fig:ex-dim1}
\end{figure}
\begin{ex}[rasped line segments] \label{ex:dim1}
Let $\mb{A} = \bR$, $\mb{S} =\{-1,1\}$, $X= [a,b] \subset \mb{A}$, and $x \in X$. 
Since the set $\frac{X-\{x\}}{r} \cap \mb{S}$ is constant for all sufficiently small $r>0$ for each $x \in X$, the limit is
\begin{equation*}
\lim_{r \searrow 0}\frac{X-\{x\}}{r} \cap \mb{S}
=
\begin{cases}
\{+1\},& x = a \\
\{-1,+1\},& a < x < b \\
\{-1\}, & x = b.
\end{cases}
\end{equation*}
Therefore the recipe \eqref{eq:recipe} gives
\begin{align} \label{eq:rasped-interval}
\sto{X} 
& \,=\, \{(a;+1)\} 
\,\sqcup\,\{(x;\pm1):a<x<b\} 
\,\sqcup\,\{(b; -1)\} \\  \nonumber
& \, \subset \stok{\mb{A}}{1} = \bR \times \{-1,1\},
\end{align}
as foreseen in Step 1 of Section \ref{sec:motiv}.

This rasped line segment is indeed divided seamlessly.
Given $c \in \bR$ such that $a < c <b$, the rasped representations of $X_1 := [a,c], X_2 := [c,b] \subset \mb{A}$ are
\begin{align*}
\sto{X}_1 
&\,=\, \{(\stackrel{\phantom{\_\_}}{a};+1)\} 
\,\sqcup\, \{(x;\pm1):~\stackrel{\phantom{\_\_}}{a}~<x<~\stackrel{\phantom{\_\_}}{c}\} 
\,\sqcup\, \{(\stackrel{\phantom{\_\_}}{c}; -1)\}, \\ 
\sto{X}_2
&\,=\, \{(\stackrel{\phantom{\_\_}}{c};+1)\} 
\,\sqcup\, \{(x;\pm1):~\stackrel{\phantom{\_\_}}{c}~<x<~\stackrel{\phantom{\_\_}}{b}\} 
\,\sqcup\, \{(\stackrel{\phantom{\_\_}}{b}; -1)\},
\end{align*}
and thus we have $\sto{X}_1 \sqcup \sto{X}_2 = \sto{X}$ (see Figure \ref{fig:ex-dim1}).

Similar results can be observed for arcs in a $1$-dimensional circle as well. 
\end{ex}
\begin{figure}
\centering
\begin{subfigure}{0.42\textwidth}
\centering
\includegraphics{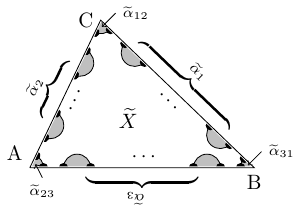}
\caption{\label{subfig:ex-dim2a}}
\end{subfigure}
\qquad
\begin{subfigure}{0.42\textwidth}
\centering
\includegraphics{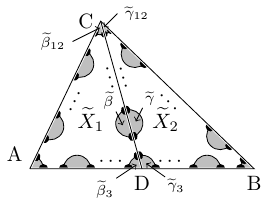}
\caption{\label{subfig:ex-dim2b}}
\end{subfigure}
\caption{A rasped triangle and its division.}
\label{fig:ex-dim2}
\end{figure}

\begin{ex}[rasped triangles] \label{ex:dim2}
Let $\mb{A} = \bR^2$, $\mb{S}=\{(x,y)\in \bR^2 : x^2 +y^2 = 1\}$, and $X$ be the closed triangle with vertices $\mr{A}$, $\mr{B}$, $\mr{C} \in \mb{A}$ in counterclockwise order. 
Let us denote by $(\mr{BC}), (\mr{CA})$ and $(\mr{AB})$ the interior of each line segment and by $(\mr{ABC})$ the interior of the triangle $\mr{ABC}$.
Again, the set $\frac{X-\{x\}}{r} \cap \mb{S}$ is constant for all sufficiently small $r>0$ for each $x \in X$, and so the limit is (see Figure \ref{subfig:ex-dim2a})
\begin{equation*}
\lim_{r \searrow 0}\frac{X-\{x\}}{r} \cap \mb{S}
= 
\begin{cases}
\alpha_{23}:=[\unitv{\mr{AB}},\unitv{\mr{AC}}], & x = \mr{A} \\
\alpha_{31}:=[\unitv{\mr{BC}},\unitv{\mr{BA}}], & x = \mr{B} \\
\alpha_{12}:=[\unitv{\mr{CA}},\unitv{\mr{CB}}], & x = \mr{C} \\
\alpha_{1\phantom{1}}:=[\unitv{\mr{BC}},-\unitv{\mr{BC}}], & x \in (\mr{BC}) \\
\alpha_{2\phantom{1}}:=[\unitv{\mr{CA}},-\unitv{\mr{CA}}], & x \in (\mr{CA}) \\
\alpha_{3\phantom{1}}:=[\unitv{\mr{AB}},-\unitv{\mr{AB}}], & x \in (\mr{AB}) \\
\mb{S}, & x \in (\mr{ABC}).
\end{cases}
\end{equation*}
Therefore the recipe \eqref{eq:recipe} gives
\begin{align*} 
\sto{X}
& \,{}={}\, \{\mr{A}\} \times \sto{\alpha}_{23} 
\,\sqcup\, \{\mr{B}\}\times \sto{\alpha}_{31} 
\,\sqcup\, \{\mr{C}\}\times \sto{\alpha}_{12} \\
& \phantom{\,{}={}\,}
\,\sqcup\, (\mr{BC}) \times \sto{\alpha}_1 
\,\sqcup\, (\mr{CA})\times \sto{\alpha}_2 
\,\sqcup\, (\mr{AB})\times\sto{\alpha}_3 
\,\sqcup\, (\mr{ABC}) \times \stok{\mb{S}}{1} \\
& \,{}\subset{}\, \stok{\mb{A}}{2} = \mb{A} \times \stok{\mb{S}}{1}.
\end{align*}

This rasped triangle is divided seamlessly into two rasped triangles as follows. 
Let $\mr{D}$ be a point on $\mr{AB}$ and ${X}_1$ and ${X}_2$ the triangles $\mr{ADC}$ and $\mr{BDC}$, respectively. 
Then we likewise have
\begin{align*}
\sto{X}_1
& \,{}={}\, \{\mr{A}\} \times \sto{\alpha}_{23} 
\,\sqcup\, \{\mr{D}\}\times \sto{\beta}_{3} 
\,\sqcup\, \{\mr{C}\}\times \sto{\beta}_{12} \\
& \phantom{\,{}={}\,}\ \sqcup\ (\mr{DC}) \times \sto{\beta} 
\,\sqcup\, (\mr{CA})\times \sto{\alpha}_2 
\,\sqcup\, (\mr{AD})\times\sto{\alpha}_3 
\,\sqcup\, (\mr{ADC}) \times \stok{\mb{S}}{1}
\\
\sto{X}_2
& \,{}={}\, \{\mr{D}\} \times \sto{\gamma}_{3} 
\,\sqcup\, \{\mr{B}\}\times \sto{\alpha}_{31} 
\,\sqcup\, \{\mr{C}\}\times \sto{\gamma}_{12} \\
& \phantom{\,{}={}\,}\ \sqcup\ (\mr{DC}) \times \sto{\gamma} 
\,\sqcup\, (\mr{BC})\times \sto{\alpha}_1 
\,\sqcup\, (\mr{DB})\times\sto{\alpha}_3 
\,\sqcup\, (\mr{DBC}) \times \stok{\mb{S}}{1}
\end{align*}
for appropriate arcs ${\beta}$, ${\beta}_{3}$, ${\beta}_{12}$, ${\gamma}$, ${\gamma}_3$, and ${\gamma}_{12}$. 
Then, the equality $\sto{X} = \sto{X}_1 \sqcup \sto{X}_2$ follows from seamless division (see Figure \ref{subfig:ex-dim2b})
\begin{equation*} \label{eq:ex-angle-ahead}
\sto{\beta}_3\sqcup\sto{\gamma}_3 = \sto{\alpha}_3, \quad
\sto{\beta}_{12} \sqcup \sto{\gamma}_{12} = \sto{\alpha}_{12},\quad
\sto{\beta}\sqcup\sto{\gamma} = \stok{\mb{S}}{1}
\end{equation*}
of arcs similar to the one in Example \ref{ex:dim1}.
\end{ex}
\begin{figure}
\centering
\includegraphics{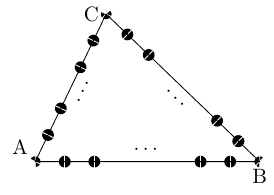}
\caption{A rasped triangular skeleton in a plane.}
\label{fig:ex-dim12}
\end{figure}
\begin{ex}[continued]\label{ex:rasp-1d2}
Let $Y$ be the skeleton of the triangle $\mr{ABC}$,
\begin{equation*}
Y = \{\mr{A}\} \sqcup \{\mr{B}\} \sqcup \{\mr{C}\} 
\sqcup (\mr{AB}) \sqcup (\mr{BC}) \sqcup (\mr{CA}).
\end{equation*}
For the same reason as in similar calculations in previous examples, we have
\begin{equation*}
\lim_{r \searrow 0}\frac{Y-\{x\}}{r} \cap \mb{S}
= 
\begin{cases}
\{\unitv{\mr{AB}},\unitv{\mr{AC}}\}, & x = \mr{A} \\
\{\unitv{\mr{BC}},\unitv{\mr{BA}}\}, & x = \mr{B} \\
\{\unitv{\mr{CA}},\unitv{\mr{CB}}\}, & x = \mr{C} \\
\{\unitv{\mr{BC}},-\unitv{\mr{BC}}\}, & x \in (\mr{BC}) \\
\{\unitv{\mr{CA}},-\unitv{\mr{CA}}\}, & x \in (\mr{CA}) \\
\{\unitv{\mr{AB}},-\unitv{\mr{AB}}\}, & x \in (\mr{AB}),
\end{cases}
\end{equation*}
and the recipe \eqref{eq:recipe} yields (see Figure \ref{fig:ex-dim12})
\begin{align*}\textstyle
\sto{Y} 
 &\,{}={}\, \{\mr{A}\} \times \{\unitv{\mr{AB}},\unitv{\mr{AC}}\}
\,\sqcup\, \{\mr{B}\} \times \{\unitv{\mr{BC}},\unitv{\mr{BA}}\} 
\,\sqcup\, \{\mr{C}\} \times \{\unitv{\mr{CA}},\unitv{\mr{CB}}\} \\ \textstyle
 &\phantom{\,{}={}\,} 
\,\sqcup\, (\mr{BC}) \times \{\unitv{\mr{BC}},-\unitv{\mr{BC}}\}
\,\sqcup\, (\mr{CA}) \times \{\unitv{\mr{CA}},-\unitv{\mr{CA}}\}
\,\sqcup\, (\mr{AB}) \times \{\unitv{\mr{AB}},-\unitv{\mr{AB}}\} \\
& \,{}\subset{}\, \stok{\mb{A}}{1} = \mb{A} \times \stok{\mb{S}}{0} = \mb{A} \times \mb{S}
\end{align*} 
Note that $\sto{Y} \subset \stok{\mb{A}}{1}$ is equal to the disjoint union of three rasped lines:
\begin{align*}\textstyle
\sto{Y} 
&\,{}={}\,
\{\mr{B}\} \times \{\unitv{\mr{BC}}\} 
\,\sqcup\, (\mr{BC}) \times \{\unitv{\mr{BC}},\unitv{\mr{CB}}\}
\,\sqcup\, \{\mr{C}\} \times \{\unitv{\mr{CB}}\}  \\ \textstyle
&\,\phantom{{}={}}\, 
\,\sqcup\, \{\mr{C}\} \times \{\unitv{\mr{CA}}\}
\,\sqcup\, (\mr{CA}) \times \{\unitv{\mr{CA}},\unitv{\mr{AC}}\}
\,\sqcup\, \{\mr{A}\} \times \{\unitv{\mr{AC}}\} \\ \textstyle
&\,\phantom{{}={}}\, 
\,\sqcup\, \{\mr{A}\} \times \{\unitv{\mr{AB}}\}
\,\sqcup\, (\mr{AB}) \times \{\unitv{\mr{AB}},\unitv{\mr{BA}}\}
\,\sqcup\, \{\mr{B}\} \times \{\unitv{\mr{BA}}\} \\ \textstyle
&{\,{}={}\,} 
\sto{\mr{BC}} \,\sqcup\, \sto{\mr{CA}} \,\sqcup\, \sto{\mr{AB}}.
\end{align*} 
\end{ex}

\begin{ex}
Wallace-Bolyai-Gerwien theorem states that any planar polygon is characterized up to scissor congruence by its area: two polygons $P$, $Q$ in $\bR^2$ have the same area if and only if there exist decompositions $P = \bigcup_{n=1}^N P_i$ and $Q = \bigcup_{n=1}^N Q_i$ into triangles such that $P_i$ is congruent to $Q_i$. 
Now, in terms of rasped polygons, we can reformulate the theorem as follows.
\begin{quote}
Two polygons $P$, $Q$ in $\bR^2$ have the same area if and only if there exist decompositions $\sto{P} = \bigsqcup_{n=1}^N \sto{P}_i$ and $\sto{Q} = \bigsqcup_{n=1}^N \sto{Q}_i$ of rasped polygons into rasped triangles such that $P_i$ is congruent to $Q_i$. 
\end{quote}
This formulation can be proved by copying original proofs (e.g., \cite[\S24]{Hartshorne2000}) with rasped polygons or by Theorem \ref{thm:rasp-curv}.
\end{ex}

We next consider rasped representation of curvilinear objects.

\begin{figure}
\centering
\begin{subfigure}{0.48\textwidth}
\centering
\includegraphics{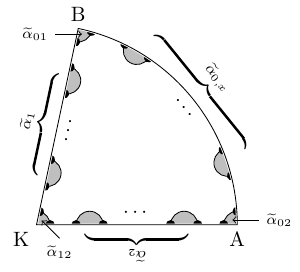}
\caption{}
\label{subfig:rasped-sector-a}
\end{subfigure}
\begin{subfigure}{0.48\textwidth}
\centering
\includegraphics{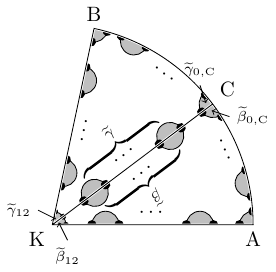}
\caption{}
\label{subfig:rasped-sector-b}
\end{subfigure}
\caption{A rasped sector in a plane.}
\label{fig:rasped-sector}
\end{figure}

\begin{ex}[rasped sectors]\label{ex:rasp-sector}
Let $\mb{A} = \bR^2$ and $X \subset \mb{A}$ be a sector corresponding to an arc $\mr{AB}$ of a circle centered at $\mr{K}$. 
If $x \in X$ but $x$ is not on the arc $\mr{AB}$, we have
\begin{equation*}
\lim_{r\searrow 0}\frac{X-\{x\}}{r} \cap \mb{S}
= 
\begin{cases}
\alpha_{12} := [\unitv{\mr{KA}},\unitv{\mr{KB}}], & x =\mr{K} \\
\alpha_{2\phantom{0}} := [\unitv{\mr{KA}},-\unitv{\mr{KA}}], & x \in (\mr{KA}) \\
\alpha_{1\phantom{0}} := [-\unitv{\mr{KB}},\unitv{\mr{KB}}], & x \in (\mr{KB}) \\
\mb{S}, & x \in (\text{sector }\mr{KAB}).
\end{cases}
\end{equation*}
(The notation $(\text{sector }\mr{KAB})$ denotes the interior of the sector $\mr{KAB}$.)
If, on the other hand, $x$ is on the arc ${\mr{AB}}$, each set $\frac{X-\{x\}}{r} \cap \mb{S}$ increases as $r \searrow 0$, and we have
\begin{equation*}
\lim_{r\searrow 0}\frac{X-\{x\}}{r} \cap \mb{S}
= 
\begin{cases}
\alpha_{02\phantom{,}} := [\unitv{\mr{KA}}^\ast,\unitv{\mr{AK}}^{}], & x =\mr{A} \\
\alpha_{0,x} := [\unitv{\mr{K}x}^\ast,\unitv{x\mr{K}}^\ast], & x \in (\text{arc } \mr{AB}) \\
\alpha_{01\phantom{,}} := [\unitv{\mr{BK}}^{},\unitv{\mr{BK}}^\ast], & x =\mr{B},
\end{cases}
\end{equation*}
(The notation $(\text{arc } \mr{AB})$ denotes the interior of the arc $\mr{AB}$.) where $\ast$ denotes the counterclockwise rotation of the given vector by the right angle. 
Here, the limits are closed intervals rather than open intervals because we consider rasped representation of closed subsets only, as mentioned at the beginning of this section. 
Then recipe \eqref{eq:recipe} yields
\begin{align*}
\sto{X} 
& \,=\, \{\mr{K}\} \times \sto{\alpha}_{12}
\,\sqcup\, \{\mr{A}\} \times \sto{\alpha}_{02}
\,\sqcup\, \{\mr{B}\} \times \sto{\alpha}_{01} 
\,\sqcup\, (\mr{KA}) \times \sto{\alpha}_{2}
\,\sqcup\, (\mr{KB}) \times \sto{\alpha}_{1} \\
&\phantom{\,{}={}\,} \,\sqcup\, \bigsqcup_{x \in (\text{arc } \mr{AB})} \{x\} \times \sto{\alpha}_{0,x}
\,\sqcup\, (\text{sector } \mr{KAB}) \times \stok{\mb{S}}{1} \\
& \,\subset\, \stok{\mb{A}}{2} = \mb{A} \times \stok{\mb{S}}{1} .
\end{align*}

This rasped sector is divided seamlessly into two sectors as follows. 
Let $\mr{C}$ be a point on the arc $\mr{AB}$ and ${X}_1$ and ${X}_2$ the sectors corresponding to arcs $\mr{AC}$ and $\mr{CB}$, respectively. 
Then we likewise have
\begin{align*}
\sto{X}_1 
& \,=\, \{\mr{K}\} \times \sto{\beta}_{12}
\,\sqcup\, \{\mr{A}\} \times \sto{\alpha}_{02}
\,\sqcup\, \{\mr{C}\} \times \sto{\beta}_{0,\mr{C}} 
\,\sqcup\, (\mr{KA}) \times \sto{\alpha}_{2}
\,\sqcup\, (\mr{KC}) \times \sto{\beta} \\
&\phantom{\,{}={}\,} \,\sqcup\, \bigsqcup_{x \in (\text{arc } \mr{AC})} \{x\} \times \sto{\alpha}_{0,x}
\,\sqcup\, (\text{sector } \mr{KAC}) \times \stok{\mb{S}}{1}\\
\sto{X}_2 
& \,=\, \{\mr{K}\} \times \sto{\gamma}_{12}
\,\sqcup\, \{\mr{C}\} \times \sto{\gamma}_{0,\mr{C}}
\,\sqcup\, \{\mr{B}\} \times \sto{\alpha}_{01} 
\,\sqcup\, (\mr{KC}) \times \sto{\gamma}
\,\sqcup\, (\mr{KB}) \times \sto{\alpha}_{1} \\
&\phantom{\,{}={}\,} \,\sqcup\, \bigsqcup_{x \in (\text{arc } \mr{CB})} \{x\} \times \sto{\alpha}_{0,x}
\,\sqcup\, (\text{sector } \mr{KCB}) \times \stok{\mb{S}}{1},
\end{align*}
for appropriate arcs $\beta$, $\beta_{0,\mr{C}}$, $\beta_{12}$, $\gamma$, $\gamma_{0,\mr{C}}$, and $\gamma_{12}$ (see Figure \ref{subfig:rasped-sector-b}). 
The equality $\sto{X} = \sto{X}_1 \sqcup \sto{X}_2$ follows from seamless division
\begin{equation*}
\sto{\beta} \sqcup \sto{\gamma} = \stok{\mb{S}}{1},
\quad \sto{\beta}_{0,\mr{C}} \sqcup \sto{\gamma}_{0,\mr{C}} = \sto{\alpha}_{0,\mr{C}},
\quad \sto{\beta}_{12} \sqcup \sto{\gamma}_{12} = \sto{\alpha}_{12}.
\end{equation*}
of rasped arcs.
\end{ex}

\begin{figure}
\centering
\includegraphics{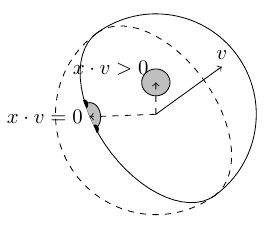}
\caption{A rasped half-sphere.}
\label{fig:rasped-half-sphere}
\end{figure}
\begin{ex}[rasped half-spheres]\label{ex:rasp-half-sphere}
Let $\mb{A} = \mb{V} = \bF^3$, $v \in \mb{V}$ a unit vector, and $X$ be a half-sphere $\mb{S} \cap \{w \in \mb{V}:w\cdot v \geq 0\} \subset \mb{V}$. 
Then for each $x \in X$ we have (see Figure \ref{fig:rasped-half-sphere}) 
\begin{equation*}
\lim_{r \searrow 0}\frac{X-\{x\}}{r} \cap \mb{S}
=
\begin{cases}
\{w \in \mb{S}: x \cdot w = 0 \} , & x\cdot v > 0, \\
\{w \in \mb{S}: x \cdot w = 0, v\cdot w \geq 0\}, &x \cdot v = 0.
\end{cases}
\end{equation*}
That is, the limit is a full circle for $x$ inside $X$ and a half-circle for $x$ on the boundary of $X$. 
The rasped representation of these circles and half-circles can be obtained similarly as in Example \ref{ex:dim1}.
Then recipe \eqref{eq:recipe} yields
\begin{align*}
\sto{X} 
&\,=\, \{(x,v_1,v_2) \in \stok{\mb{S}}{2} : x \cdot v > 0 \} \\
& \phantom{\,=\,} \,\sqcup\,\{ (x, v_1, v_2) \in \stok{\mb{S}}{2} : x \cdot v = 0 \text{ and } \big( v_1 \cdot v > 0 \text{ or } (v_1 \cdot v=0, v_2 \cdot v > 0)\big)\}
\end{align*}
Similarly, the rasped representation of the complementary half sphere $Y = -X$ is
\begin{align*}
\sto{Y}
&\,=\, \{(x,v_1,v_2) \in \stok{\mb{S}}{2} : x \cdot v < 0 \} \\
& \phantom{\,=\,} \,\sqcup\,\{ (x, v_1, v_2) \in \stok{\mb{S}}{2} : x \cdot v = 0 \text{ and } \big( v_1 \cdot v < 0 \text{ or } (v_1 \cdot v=0, v_2 \cdot v < 0)\big)\},
\end{align*}
and we observe seamless division: $\sto{X} \sqcup \sto{Y} = \stok{\mb{S}}{2} \subset \stok{\mb{V}}{2}$. 
\end{ex}

\begin{figure}
\centering
\begin{subfigure}{0.325\textwidth}
\centering
\includegraphics{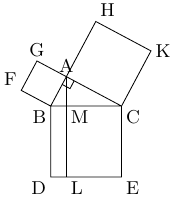}
\caption{}
\label{subfig:elem-1-47}
\end{subfigure}
\begin{subfigure}{0.325\textwidth}
\centering
\includegraphics{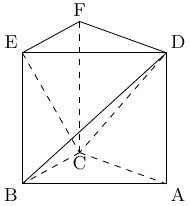}
\caption{}
\label{subfig:elem-12-7}
\end{subfigure}
\begin{subfigure}{0.325\textwidth}
\centering
\includegraphics{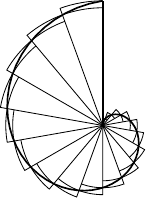}
\caption{}
\label{subfig:archim-spiral}
\end{subfigure}
\caption{Diagrams for selected classical Greek geometrical theorems. 
(a) Pythagorean theorem (\textit{Elements} I.47) 
(b) \textit{Elements} XII.7
(c) Area bounded by Archimedean spiral (\textit{On Spirals} 21, 24)
}
\label{fig:ex-real-greek}
\end{figure}

\paragraph{Rasped representation for classical Greek geometry}{
Rasped representation fits various situations in classical Greek geometry texts. 

In the proof of the Pythagorean theorem in \textit{Elements} I.47 the square $\mr{BDEC}$ is equal to the sum of the rectangle $\mr{BDLM}$ and the rectangle $\mr{CELM}$ (see Figure \ref{subfig:elem-1-47}). 
And \textit{Elements} XII.7 explains that the triangular prism $\mr{ABC}\mr{DEF}$ is divided into three triangular pyramids $\mr{CABD}$, $\mr{DBCE}$, and $\mr{DCEF}$ having equal volumes (see Figure \ref{subfig:elem-12-7}).\footnote{%
The Greek text states that the three triangular pyramids are equal. 
See note \ref{fn:equal} for the specific meaning of the word ``equal'' in \textit{Elements}.
} 
Indeed, as in the above example, a rasped square can indeed be the sum (disjoint union) of two rasped rectangles, and a rasped prism can be divided into (disjoint union of) rasped pyramids. 

Rasped representation is also suitable for Archimedes' approximation of the region bounded by a spiral in \textit{On Spirals} (see Figure \ref{subfig:archim-spiral}). 
How could the spiral region be greater than the sum of inner sectors and less than the sum of the outer sectors, without replacing those objects by their areas, if these are described as closed or open subsets of $\bR^2$?
}

\paragraph{Limitations?}{%
The above examples have already shown some features of the rasped representation compared with the conventional one. 
First, it is uncertain whether, or probably not, every subset of $\bR^d$ can be applied to rasped representation, which means that rasped representation can describe a smaller variety of geometric objects. 
Second, objects of different dimensions correspond, under rasped representation, to subsets of different ambient spaces. 
For example, (rasped) points, rasped lines, rasped polygons in $\mb{A}$ are subsets of $\stok{\mb{A}}{0}(=\mb{A})$, $\stok{\mb{A}}{1}$, and $\stok{\mb{A}}{2}$, respectively. 
Accordingly, some definitions should be modified to be more complicated than usual: the intersection point of two rasped lines $\sto{X}$ and $\sto{Y}$ is not $\sto{X} \cap \sto{Y}$, and a rasped point $\sto{X}$ can be on a rasped line $\sto{Y}$ and $\sto{X} \cap \sto{Y} = \varnothing$ at the same time. 

As will be discussed in Section \ref{sec:rasp-curv}, rasped representation applies to a fairly large range of subsets of $\bR^d$ including many familiar geometric objects, and no two different subsets share the same representation (i.e., rasped representation is injective). 
Thus, many geometric theories can be described by rasped representation.
}

\section{Definitions and proofs}\label{sec:rasp-curv}

In this section, we clarify the precise meaning of the rasped representation and the recipe \eqref{eq:recipe} and prove seamless division.
We begin with the definitions of the range $\stok{\mb{M}}{k}$ and the domain $\fS_d(\mb{M})$ of rasped representation, for a Riemannian manifold $\mb{M} = (\mb{M},\langle \cdot,\cdot\rangle)$.

\begin{figure}
\centering
\includegraphics{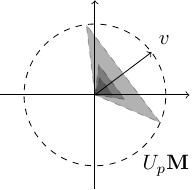}
\caption{Shape of $\gdc_{v,\delta}$ for $\delta=0.5$, $0.3$, and $0.1$.}
\label{fig:gdc}
\end{figure}

Given $x \in M$, let
\begin{equation*}
U_x \mb{M} := \{ v \in T_x \mb{M}: \langle v,v \rangle  = 1\}
\end{equation*}
be the set of unit tangent vectors.
In addition, given $v \in U_x\mb{M}$ and $\delta > 0$, let (see Figure \ref{fig:gdc})
\begin{equation*}
\gdc_{v,\delta}
:= \{ u \in T_x\mb{M}: (1-\delta) \|u\| \|v\| < u\cdot v < \delta \}.
\end{equation*}
\begin{defn}\label{def:sto-mani}
Let $\mb{M}=(\mb{M},\langle \cdot,\cdot\rangle)$ be a Riemannian $d$-manifold and $k$ an integer such that $0 \leq k \leq d$.
\begin{enumerate}
\item We define $\stok{\mb{M}}{k}$ as
\begin{align*}
\stok{\mb{M}}{k} 
& = 
\begin{cases}
\mb{M}, & k = 0 \\
\bigsqcup_{x \in \mb{M}} \{x \} \times \stok{U_x \mb{M}}{k-1}, & k > 0
\end{cases}
\end{align*}
In other words, for $k>0$ we have
\begin{align*}
\stok{\mb{M}}{k}=\{(x;{}&{}v_1,\cdots,v_k) \in \mb{M} \times (T\mb{M})^k : \\ &v_1,\cdots,v_k \in T_x\mb{M} \text{ are orthonormal}\},
\end{align*}

\item \label{item:sto-mani}
We define $\mathfrak{S}_d(\mb{M})$ as the set of all subsets $X\subset \mb{M}$ satisfying the following:
\begin{enumerate}
\item $X = \mathbf{cl}(\mathbf{int}\, X)$.

\item For each $x \in X$, the subset
\begin{equation*}
\mathcal{C}_{X,x}:=
\{v \in U_x \mb{M}: \exists \delta > 0
\,( \exp_x(\gdc_{v,\delta}) \subset X 
\text{ or } \exp_x(\gdc_{v,\delta}) \subset \mb{M} - X )\}
\end{equation*}
has dense interior in $U_x M$. 

\item For each $x \in X$, the set
\begin{equation*}
U_x X := \{v \in U_x \mb{M} : v = \gamma'(0) \text{ for a smooth curve }\gamma:[0,\varepsilon_0) \rightarrow X \}
\end{equation*}
of unit tangent vectors of $X$ at $x$ satisfies $\ci(U_x X) \in \mathfrak{S}_{d-1}(U_x \mb{M})$.
\end{enumerate}
\end{enumerate}
\end{defn}
The set $\stok{\mb{M}}{k}$ was defined similarly to Definition \ref{def:sto-basic}. 
To figure out which subsets of $\mb{M}$ belong to $\fS_d(\mb{M})$, we need the following two lemmas, which will be proved later.

\begin{lem}\label{lem:sublvl-set}
If $f:\mb{M} \rightarrow \bR$ is a smooth map having $0$ as a regular value, then $X := \{p \in \mb{M}: f(p) \leq 0\} \in \fS_d({\mb{M}})$. 
\end{lem}

\begin{lem} \label{lem:curv-union-inter}
If $X,Y \in \fS_d({M})$, then $X \cup Y$,$\ci(X \cap Y)\in \fS_d({M})$.
\end{lem}

By applying Lemma \ref{lem:sublvl-set} and \ref{lem:curv-union-inter} together (iteratively if necessary), we can prove that polygons, sectors, areas bounded by conic sections and straight lines belong to $\fS(\bR^2)$, and convex polytopes, balls, cylinders belong to $\fS_d(\bR^3)$. 

We now state our main theorem on rasped representation and seamless division.

\begin{thm} \label{thm:rasp-curv}
Let $\mb{M}$ be a Riemannian $d$-manifold. 
The recursive formula
\begin{equation} \label{eq:rasp-curv}
X \,\mapsto\, \sto{X} := 
\begin{cases}
X , & k = 0\\
\bigsqcup_{x \in X} \{x\} \times \sto{\ci(U_x X)}, & k > 0.
\end{cases}
\end{equation}
defines an injective map on $\fS_d(M)$ whose image consists of subsets of $\stok{\mb{M}}{d}$ and satisfying the following properties.
\begin{enumerate}
\item $\sto{X} \cup \sto{Y} = \sto{X \cup Y}$.
\item $\sto{X} \cap \sto{Y} = \sto{\ci(X \cap Y)}$.
\end{enumerate}
\end{thm}

\begin{cor}[seamless division]
If $X$,$Y \in \fS_d(\mb{M})$ and $\interior (X \cap Y) = \varnothing$, then we have $\sto{X \cup Y} = \sto{X} \sqcup \sto{Y}$.
\end{cor}

\begin{ex}
If $f:\mb{M} \rightarrow \bR$ is a smooth map having $0$ as a regular value, its rasped representation $\sto{X}$ is
\begin{equation*}
\begin{aligned}
\sto{X} = \big\{& (x;v_1,\cdots,v_d) \in \stok{\mb{M}}{d} :
\\
&\big(f(x),df_p(v_1),\cdots,df_p(v_1+\cdots,v_d)\big) \prec (0,0,\cdots,0) \big\},
\end{aligned}
\end{equation*}
where $\prec$ denotes the lexicographical comparison: For two tuples of real numbers of the same length, $(a_1,\cdots,a_k) \prec (b_1,\cdots,b_k)$ if and only if there exists an index $i_0$ such that $a_i = b_i$ for $i<i_0$ and $a_{i_0} < b_{i_0}$. 
The above expression can be deduced by mathematical induction. 
Moreover, Examples \ref{ex:dim1}--\ref{ex:rasp-half-sphere} are reproduced by taking suitable intersections of such spaces and applying Theorem \ref{thm:rasp-curv}. 
\end{ex}

We now turn to the proofs of the lemmas and the main theorem. 

\begin{lem} \label{lem:ci-party}
Given $X,Y \in \fS_d({M})$ and $x \in X \cap Y$, we have
\begin{enumerate}
\item $\ci(U_x X) \cup \ci(U_x Y) = \ci (U_x (X \cup Y))$
\item $\ci(\ci(U_x X) \cap \ci(U_x Y)) = \ci (U_x (\ci(X \cap Y)))$
\end{enumerate}
\end{lem}
\begin{proof}
\begin{enumerate}
\item The ``$\subset$'' is implied from
\begin{align*}
U_x X, U_x Y \subset U_x (X \cup Y) 
& \,\Rightarrow\, \interior(U_x X), \interior(U_x Y) \subset \interior (U_x (X \cup Y)) \\
& \,\Rightarrow\, \interior(U_x X) \cup \interior(U_x Y) \subset \interior (U_x (X \cup Y))
\end{align*}
and
\begin{align*}
\ci(U_x X) \cup \ci (U_x Y) 
& = \closure( \interior(U_x X) \cup \interior(U_x Y))\\
& \subset \closure (\interior(U_x (X\cup Y))).
\end{align*}

For the ``$\supset$'' part, we have
\begin{align*}
\closure(\interior(U_x(X \cup Y)))
& = \closure (\interior (U_x (X \cup Y)) \cap \mathcal{C}_{X,x} \cap \mathcal{C}_{Y,x} ).
\end{align*}
because $\mathcal{C}_{X,x} \cap \mathcal{C}_{Y,x}$ contains an open dense subset of $T_x \mb{M}$.\footnote{%
Here we apply an easily proved lemma: If $A$ is an open subset and $B$ is an open dense subset of a topological space, then $\closure(A \cap B) = \closure(A)$.
} 
Once we prove
\begin{equation*}
\interior (U_x (X \cup Y)) \cap \mathcal{C}_{X,x} \cap \mathcal{C}_{Y,x} 
\,\subset\, \interior(U_x (X \cup Y))
\end{equation*}
the proof of ``$\supset$'' part will be complete. 
If $v$ belongs to the left-hand side, then $v \in \interior (U_x (X\cup Y))$ implies that two conditions $\exp_x(\gdc_{v,\delta}) \subset \mb{M}-X$ and $\exp_x(\gdc_{v,\delta}) \subset \mb{M}-Y$ cannot hold simultaneously. 
Then, $v \in \mathcal{C}_{X,x} \cap \mathcal{C}_{Y,x}$ implies that $\exp_x(\gdc_{v,\delta}) \subset X$ or $\exp_p(\gdc_{v,\delta}) \subset Y$, which in turn gives $v \in \interior(U_x (X \cup Y))$. 

\item The ``$\supset$'' part is implied from
\begin{align*}
U_x X, U_x Y \supset U_x \ci (X \cap Y)
&\,\Rightarrow\, \interior(U_x X), \interior(U_x Y) \supset \interior(U_x \ci (X \cap Y)) \\
&\,\Rightarrow\, \interior(U_x X) \cap \interior(U_x Y) \supset \interior(U_x \ci (X \cap Y))
\end{align*}
and
\begin{align*}
\interior(\ci(U_x X) \cap \ci(U_x Y))
& = \ic(\interior(U_x X)) \cap \ic(\interior(U_x Y)) \\
& \supset \interior(U_x X) \cap \interior(U_x Y) \\
& \supset \interior(U_x \ci(X \cap Y)).
\end{align*}
For the ``$\subset$'' part, we have
\begin{align*}
\ci(\ci(U_x X) \cap \ci(U_x Y)) 
& = \closure(\interior(\ci(U_x X) \cap \ci(U_x Y)) \cap \mathcal{C}_{X,x} \cap \mathcal{C}_{Y,x}) \\
& = \closure(\interior(\ci(U_x X)) \cap \interior(\ci(U_x Y)) \cap \mathcal{C}_{X,x} \cap \mathcal{C}_{Y,x}) 
\end{align*}
because $\mathcal{C}_{X,x} \cap \mathcal{C}_{Y,x}$ contains an open dense subset of $T_x \mb{M}$. 
Once we prove
\begin{equation*}
\interior(\ci(U_x X)) \cap \interior(\ci(U_x Y)) \cap \mathcal{C}_{X,x} \cap \mathcal{C}_{Y,x} 
\,\subset\, \interior(U_x\ci(X \cap Y))
\end{equation*}
the proof of the ``$\subset$'' part will be complete. 
If $v$ belongs to the left-hand side, then $v \in \interior(\ci(U_x X)) \cap \interior(\ci(U_x Y))$ implies that neither $\exp_x (\gdc_{v,\delta}) \subset \mb{M}-X$ nor $\exp_x (\gdc_{v,\delta}) \subset \mb{M}-Y$ can hold. 
Then, $v \in \mathcal{C}_{X,x} \cap \mathcal{C}_{Y,x}$ implies that $\exp_x (\gdc_{v,\delta}) \subset X \cap Y$, which in turn gives $v \in \interior(U_x\ci(X \cap Y))$. 
\end{enumerate}
\end{proof}

\begin{proof}[Proof of Lemma \ref{lem:sublvl-set}]
We check each of conditions (a)--(c) of Definition \ref{def:sto-mani}-\ref{item:sto-mani}.
\begin{enumerate}
\item[(a)] $X \supset \closure(\interior\,X)$ is immediate. 
$X \subset \closure(\interior\,X)$ holds because every $x \in \mb{M}$ such that $f(x) = 0$ is the limit of a sequence $(x_n)$ in $\mb{M}$ such that $f(x_n) < 0$ by nonsingularity of $f$ at $x$. 
\item[(b)] We have $T_p \mb{M} - \ker df_x \subset \mathcal{C}_{X,p}$. 
\item[(c)] $U_x X = \{df_x \leq 0 \}$ is the half-sphere of $U_x \mb{M}$, which belongs to $\fS_{d-1}(U_x \mb{M})$. 
\end{enumerate}
\end{proof}

\begin{proof}[Proof of Lemma \ref{lem:curv-union-inter}]
We check each of conditions (a)--(c) of Definition \ref{def:sto-mani}-\ref{item:sto-mani}.
\begin{enumerate}
\item[(a)] The equality $X \cup Y = \ci(X \cup Y)$ follows from
\begin{align*}
X \cup Y 
&= \closure(\interior(X)) \cup \closure(\interior(Y))\\
&= \closure(\interior(X) \cup \interior(Y)) \\
& \subset \closure(\interior(X \cup Y))
\end{align*}
and putting $K = X \cup Y$ in $K = \closure(K) \supset \closure(\interior(K))$ which holds for any closed subset $K$. 
Next, the equality $\ci(X \cap Y) = \ci(\ci(X\cap Y))$ is actually a special case of $\ci(K) = \ci(\ci(K))$ which holds for any closed subset $K$. 
This is implied by
\begin{align*}
K \supset \ci(K) 
&\,\Rightarrow\, \interior(K) \supset \interior(\ci(K))\\
&\,\Rightarrow\, \ci(K) \supset \ci(\ci(K))
\end{align*}
and
\begin{align*}
\interior(K) \subset \ci(K)
&\,\Rightarrow\, \interior(K) \subset \interior(\ci(K))\\
&\,\Rightarrow\, \ci(K) \subset \ci(\ci(K)).
\end{align*}

\item[(b)] This follows from $\mathcal{C}_{X \cup Y,x}$, $\mathcal{C}_{\ci(X \cap Y),x} \supset \mathcal{C}_{X,x} \cap \mathcal{C}_{Y,x}$

\item[(c)] This follows from Lemma \ref{lem:ci-party} and the mathematical induction on $d$. 
\end{enumerate}
\end{proof}

\begin{proof}[Proof of Theorem \ref{thm:rasp-curv}]
Since $U_x X = U_p \mb{M} \neq \varnothing$ for all $x$ in a dense subset $\interior\, X$ of $X$, we have $X = \closure(\mathrm{pr}_1(\sto{X}))$ where $\mathrm{pr}_1:\sto{M} \rightarrow M$ is the projection onto the first factor. 
This proves the injectivity. 
The remaining assertions follow from mathematical induction, recursive formula \eqref{eq:rasp-curv}, and Lemma \ref{lem:ci-party}.
\end{proof}

\section{Closing remarks}


The question at the beginning sounds quite simple, but conventional mathematical formulations had shortcomings, and a new formulation was developed to yield a more satisfactory answer. 
The new representation of geometric objects is finer than sets of points in conventional coordinate geometry, and this enabled dividing a geometric object without problems along the cut boundary. 
Because of this distinguished property, classical Greek geometry can be reinterpreted in a more natural way through the new representation. 
Especially, the absence of quantitative concepts such as length, area, and volume in classical texts is preserved even under modern set-theoretic formalization. 

The new idea of the current article, ``rasped representation'', is just another expression of elementary geometry, yet may be used for other impressive reinterpretations of known facts. 
For example, by virtue of the seamless division property, wouldn't there be interesting examples related to valuations of convex polytopes? 
This direction may be pursued in subsequent studies.

\paragraph{Acknowledgements.}{%
The author would like to thank Eunsoo Lee and Cheongmyeong Lee for their valuable advice. 
This research was supported by the Basic Science Research Program through the National Research Foundation of Korea (NRF), funded by the Ministry of Education (grant no. RS-2023-00248700). 
}

\paragraph{Conflict of interest}
The author declares no conflicts of interest.

\bibliographystyle{abbrv}
\bibliography{references.bib}
\end{document}